\documentclass{article}

\usepackage{stmaryrd,amssymb,amsmath,latexsym,amsthm,url}
\usepackage{hyperref}
\usepackage{xfrac} 
\usepackage[english]{babel}
\usepackage{eurosym}

\usepackage[usenames]{color}
\usepackage{amsthm}

\newtheorem{definition}{Definition}[section]
\newtheorem{proposition}{Proposition}[section]

\newtheorem{mydef}{Definition}

\newcommand{\quot}[2]{\ensuremath{\textstyle{{#1}/_{\!q}{#2}}}}

\newcommand{\AQ}{\mathsf{AQ}}

\newcommand{\flatfracterm}{\mathsf{flat\_fracterm}}

\newcommand{\Nat}{{\mathbb N}}
\newcommand{\Int}{{\mathbb Z}}
\newcommand{\Rat}{{\mathbb Q}}

\title{Meadow based Fracterm Theory}

\author{%
Jan A.\ Bergstra\\
{  Informatics Institute, University of Amsterdam}%
\thanks{email: \texttt{j.a.bergstra@uva.nl, janaldertb@gmail.com}.}
\date{June 6th 2019}
}

\begin{document}

\maketitle

\begin{abstract}
\noindent Fracterms are introduced as a proxy for fractions. A precise definition of fracterms is formulated and on that basis reasonably precise definitions of various classes of fracterms are given. In the context of the meadow of rational numbers viewing fractions as fracterms provides an adequate theory of fractions. A very different view on fractions is that fractions are values, i.e. rational numbers.
Fracterms are used to provide a range of intermediate definitions between these two definitions of fractions.\footnote{%
Second version of~\url{http://arxiv.org/abs/1508.01334v1}.}

\end{abstract}

\tableofcontents
\section{Introduction}
The original intention of this paper was to provide a systematic and comprehensive theory of fractions. 
However, thousands of papers have been written about fractions, admittedly mostly in an educational context, while often including significant conceptual observations, and I was unable to find a convincing structure in the 
utterly diverse literature on fractions. Although different 
perspectives on fractions have been put forward by different authors I have not seen a single paper which provides contrasting 
views on fractions in order to motivate a particular chosen perspective. 
All work concerning seems to apply a fraction talk which is perceived by the author(s) as 
self-evident and which for that reason is not in need of any elaborate justification, while making mention of the existence of
other views on fractions is not even exceptional, it seems to be running against implicit methodological principles so to say. 
Most papers explicitly discussing the concept and definition of fractions make no mention of existing work on that issue.

The question 
``how to define fractions'' is not  on the agenda, and it seems not  to be considered a well-posed question, 
being as a question only marginally less problematic than the question whether or not $1/0 = 2/0$, 
a question which many mathematicians wil reject as being ill-posed, and 
which, for that very reason is not even entitled to be answered by way of a convincing response, 
let alone a response which takes the existence of different views on the matter into account.
 
I failed to find a convincing definition of a fraction as a mathematical concept. 
This leads me to the hypothesis that in spite of its abundant use as a technical term, fraction is not a mathematical notion, a hypothesis which has been put forward in~\cite{FandinoP2007}. 
Considerable variation in the concept of fraction used by various authors, exits. 
This variation is no surprise in the light of the relativism of\cite{Shapiro2014} but nevertheless leaves room for further research: why is it the case that fraction just like proof, definition, theorem, and result, 
need not be given a rigorous definition in the presentation of arithmetic. In other words: what makes fraction 
as a notion different from say the following notions: integer, prime, rational numbers, real number, factor, field, metric space, or topological space, 
all notions for which giving rigorous definitions is standard practice. 

I wil carry out the thought experiment that fractions are defined in some meticulous detail 
and sort out the consequences of such an approach. Rather than to propose a definition 
of fractions and analyse its pro's and con's I will carry out the thought experiment that fraction is 
replaced by a precisely defined concept, in such a manner that no tension with existing terminology 
and conventions about fractions arises. Moreover, I will focus on a single proxy of fraction which comes about when 
viewing a fraction as an expression.

In order to develop a systematic way forward  I introduce ``fracterm'' as a new technical notion, where a fracterm is an expression over
an appropriate signature with division as its leading function symbol. In the absence of
pre-existing theory of fracterms an explanation of fracterms and of terminology related to fracterms can be given without risking any mismatch with existing conventions. The theory of fracterms is developed in such a manner that upon replacing fracterm 
by fraction, i.e. reading fraction as fracterm while forgetting the keyword ``fracterm'' which is then become redundant, 
an adequate theory (explanation) of fractions emerges. 
The resulting perspective on fractions, which I will refer to as ``fractions as terms'' (or fractions as fracterms)
is my preferred view on the matter. 

A plurality of different explanations of fractions can be found by reading fraction as a 
particular subclass of fracterms. Starting out with fracterms 
 helps to find some structure in the multitude  of published and unpublished views on fractions.
 \subsection{Fracterm theory, where are  is the difficulties?}
 Developing a theory of fracterms is a challenge for several reasons. 
 In giving fracterms an important role clarity is achieved by separating syntax and semantics. 
 But, doing so at every level leads to an infinite regress. 
 Thus at some level of abstraction, the desire to provide clarity by way of separating syntax 
 and semantics must give way for the desire to avoid an ongoing regress asking for more and more, increasingly futile detail.
 Below this is done by making use of the notions of explicit and implicit assimilation and explicit and implicit 
 dis-assimilation as put forward in~\cite{NicaudBG2001}. 
 
 Fracterms constitute an appearance of syntax within elementary mathematics, which might but need not replace the notion of a fraction.
 Even if fracterms do not advance to the level of serving as a replacement of fractions certain advantages of working with fracterms
 justify these investigations: (i) the notion of fracterms and fracterm terminology can be easily adapted/generalised to 
 different arithmetical datatypes such as wheels (\cite{Setzer1997,Carlstroem2004,Carlstroem2005}) transrationals (\cite{AndersonVA2002}), and transreals
 (\cite{AndersonVA2007}), 
 (ii) the perspective of fractions as fracterms is helpful when developing a theory of fractions for a specific arithmetical datatype, which matters because it is at the level of fraction theory where the different arithmetical datatypes diverge, 
 (iii) if fracterm theory turns out to be  problematic and if its development is fraught with unresolved complications, then that very observation is to some extent explanatory for how fractions are dealt with in practice.
 
\subsection{The signature of meadows}
The signature of rings consists of a sort of values, constants zero and one, and functions addition, opposite, and multiplication, with subtraction as a derived function (given by $x-y = x + (-y)$). I will assume that fixed notations come with the signature but in a more sophisticated presentation 
the symbols used for the constants and functions constitute an independent parameter. 
I will use $\_ + \_$ for addition, $-\_$ for opposite
and $\_ \cdot \_$ for multiplication, and $V$ for the sort of values.

The signature of rings is the same as the signature of fields and fields are a subclass of rings. The signature of meadows extends the signature of rings with either an inverse function (inversive meadows, meadows with inversive notation) or a division function 
(divisive meadows, meadows with divisive notation), and in both cases the other operation 
may be considered a derived operation, where
$x/y = x \cdot y^{-1}$ defines division in terms of inverse and multiplication and $x^{-1} = 1/x$ defines inverse in terms of division.
We will denote inverse with $\_^{-1}$ and division with $\_/\_$ the latter symbol optionally rendered as a horizontal bar as well.

The presence of a signature does not imply a formalistic style where syntax and semantics are distinguished in each and every notation.
On the contrary, say the object ``zero'' and its name $0$ as part of the signature may be identified, a phenomenon called assimilation 
in~\cite{NicaudBG2001}, and may be dis-assimilated at some stage when a higher conceptual resolution is needed. In the case of the constant zero, dis-assimilation leads to a ramification of notation, for instance one may keep the notation $0$ for the mathematical entity and introduce an other notation, for instance $\underline{0}$, for the name of the semantic entity in the signature. Alternatively one may insist that $0$ is a name in a signature and introduce another notation, for instance $\hat{0}$ for a semantic entity serving as its interpretation.
or conversely. The virtue of the concept of assimilation is that it explains how expressions may be used as constituting a part of the mathematical discourse, instead of merely being delegated to a parallel universe of formal utterances and sentences.

\subsection{A signature, why does it matter?}
By stipulating the existence of a signature at the same time the presence of a volume of syntax is claimed. 
For instance upon assuming the availability of the signature of divisive meadows, or upon participating in an exchange of views in
which the presence of the signature of meadows is presumed one must accept the existence of  the expressions which can be 
built from that syntax. Accepting $(1+1) \cdot (1 + 1)$ is not controversial in contrast with $1/0$ 
and $(1+1)/(1-1)$. By assuming the presence of the signature of meadows one does subscribe to the assumption that $1/0$ 
signifies some abstract entity. However, one accepts that question about the ontology of the expression $1/0$, 
cannot and must not be avoided.

Besides expressions in the signature of meadows I will make ample use of decimal notation for natural numbers, which may be 
thought of a being defined by an infinite collection of axioms $2 = 1+1$, $3 = 2+1$, etc.

\section{Rational numbers, choice of a model}
\label{Rationals}
I will assume the presence of a fixed structure for the field $\Rat$ of rational numbers. I assume that a fixed but arbitrary 
instance of a structure $\Int$ for integers is given, and all integers in the sequel of this paper  are supposed to be taken from its domain.
Natural numbers are the non-negative elements in $|\Int|$.

The domain $|\Rat|$ of $\Rat$ consists of pairs $(n,m)$ with $n,m \in \Int$ such that $ m> 0$, and $ \gcd(n,m) = 1$.
In this notation $0$is identified with $(0,1)$ and $i$ is identified with $(1,0)$. The operation $\sigma$ for simplification works as follows
on pairs $(n,m)$ with $n,m \in \Int$ and $m>0$: if $\gcd(n,m) = 1$ then $S(n,m) = (n,m)$, if for some $a,b,c \in \Int$, $n = a \cdot b$, 
$m= a \cdot c$ with $a>1,c>1$, then $\sigma(n,m) = \sigma(b,c)$. On $|\Rat|$ addition, opposite, and multiplication are given by:
$(n,m) + (a,b) = \sigma(n \cdot b+ a \cdot m, m \cdot b)$, 
$-(n,m) = (-n,m)$, and $(n,m) \cdot (a,b) = \sigma(n \cdot a, m \cdot b)$.

I consider the question ``what is a fraction'', to be independent of the question ``what is a rational number''. I do not claim that the definition of fractions requires a definitions of rational numbers as a prerequisite. The claim is instead that a discussion of fractions can be conducted 
once it has been determined what rational numbers are, and no freedom of (re)definition of rational numbers is needed for the understanding of any plausible notion of fraction. For the remainder of the paper it matters that some structure of rational numbers is chosen so that it is known what a rational number is. Clearly the notion of a fraction is not needed for defining rational numbers. We hold that fractions come in to play when working with rational numbers. Our discussion is limited to the rational numbers in the sense that I do not discuss the notion of fraction in the context of larger number structures such as 
$\mathbb{R}$ and $\mathbb{C}$. So as far as this paper is about the question ``what is a fraction?'' that question may be read as: 
given the rational numbers $\Rat$ as known: what is a fraction?

\subsection{Naturals and integers}
Natural numbers are supposed to be given and known. With the rational numbers given as above one may say that, say $7$ is an expressions which denotes $(7,1)$ and so on. I will avoid such puzzling considerations as follows. $0,1,2,3,4,\ldots,10,11, \ldots$ etc.
are numerals for natural numbers. I will assume that for natural numbers assimilation is applied to the 
name of the number and the number itself. Thus, not only $23$ denotes a natural number, it is a natural number. 
$23$ is a name, a preferred numeral, a preferred notation, and a number at the same time. 
I will deal in the same way with integers. Admittedly $-25$ denotes 
a negative integer, having applied explicit assimilation of names and meanings of integers. 

\subsection{Assimilation ground level: a matter of choice}
From a conceptual perspective assimilation of $7$ and ``the natural number denoted by $7$'' is wrong because, $7$ has a notational bias which the number itself has not. But at the same time each exposition of arithmetic will at some stage need to identify signs with their meaning. Explicit assimilation indicates that the level at which this is done comes about from an explicit mental decision, a choice which might be made in a different manner.  With the assimilation ground level I will indicate those concepts and entities for which no distinction between sign and concept is made. For the the theme of this paper, fractions, I assume that natural numbers, integers, and the well-know notations thereof are assimilated in the conventional manner. I will also assume that, say $007$ is the same as $7$, that is leading zeros for natural numbers are ignored, or skipped, automatically without further explanation. Clearly the choice of an assimilation ground 
level is a parameter, and by choosing a lower level a higher conceptual resolution may be obtained. I will also assume that assimilation 
extends to the rationals as follows: say $27$ is assimilated with $(27,1)$, and $-27$ is assimilated with $(-27,1)$. 

Following~\cite{NicaudBG2001} explicit dis-assimilation can be applied if one prefers a higher conceptual resolution. 
In the case of naturals and 
leading zeros one may introduce decimal number notations of various kinds and distinguish the numbers for 
the corresponding notations. This may be done in different ways. Taking assimilation and dis-assimilation on 
board in the methodology of 
fraction theory may seem an overkill, but I believe that it solves a fundamental difficulty. 
For an explanation of fractions, one's understanding of naturals and integers act as a sort of parameter, 
and it is difficult, if not impossible, to achieve generality with respect to 
that parameter. By adopting a deliberate policy of assimilation, i.e. of the choice of an assimilation ground level, 
the variation in that parameter may be reduced to a manageable level.

Fractions as a topic arises once it is agreed that assimilation between, 
say $2/3$ and the number denoted by that sign, is not  to be taken for granted.

\section{Fracterms and  fracterm terminology}
Central to the discussion below is the notion of a fracterm a special kind of expression.\footnote{%
In the Dutch translations, fracterm may be translated with ``breukvorm'' or ``breukterm''.}

\subsection{Arithmetical signatures}
Each notion of fracterm requires the presence of a suitable signature.
\begin{mydef}
\label{ASwD}
A an arithmetical signature  is a signature which equals or extends the signature of inversive meadows.
\end{mydef}

\begin{mydef}
\label{ASwD}
A an arithmetical signature with division is a signature which equals or extends the signature of divisive meadows.
\end{mydef}
The following definitions work under the assumption that an arithmetical signature with division $\Sigma$ is given. Precisely this
assumption is uncommon for ordinary mathematical work.
\begin{mydef}
\subsection{Primary fracterm terminology}
A definition of fracterms can be given once a suitable signature is known.
\label{Fracterm}
A fracterm is a $\Sigma$-expression with division as its leading function symbol.
\end{mydef}

The notion of a fracterm is implicitly extended with other notations for the division operator: $\frac{P}{Q}$, $ \sfrac{P}{Q}$, and $P\div Q$
each of which may also feature as the leading function symbol of a fraction.

\begin{mydef} For a fracterm $P/Q$, the expression $P$ is called the numerator and the expression $Q$ is called the denominator.
\end{mydef}
I define the quotient of a fracterm as its value. A quotient is a value, i.e. a number, and for that reason a quotient 
has no numerator and no denominator. 
\begin{mydef} The value in an arithmetical datatype of the fracterm $P/Q$ is called  the quotient of the values of $P$ and $Q$.
\end{mydef}

Moreover an integral domain may be embedded in 
a field of quotients rather than in a field of fractions.
\begin{mydef} The process of finding  the quotient of a fracterm is called dividing.
\end{mydef}

The equality sign expresses that both sides have the same value. Thus $P/Q = R/S$ expresses that the fracterms $P/Q$ and $R/S$ have the same quotient.
Every value is a quotient, given that for all terms $P$, $P/1 = P$. Consider the fracterms $P \equiv 2/3$ and $Q \equiv 4/6$, then $P$ and $Q$ are different fracterms (written $P \not \equiv Q$) and at the same time $P=Q$, i.e. $P$ and $Q$ have the same value.

The difference of two fracterms $P/Q$ and $R/S$ is the value of the subtraction $P/Q-R/S$. It follows that the difference of different fracterms can be zero (a shorthand of an be equal to zero), for instance $4/6-2/3 = 0$.

\begin{mydef} Sameness and equality:

\begin{itemize}
\item Two fracterms $P/Q$ and $R/S$ are the same, written $P/Q \equiv R/S$ if $P \equiv R$ and $Q \equiv S$, that is if respective numerators and denominators are the same terms.
 
\item Two fracterms $P/Q$ and $R/S$ are equal, written $P/Q =R/S$ if $P/Q$ and $R/ S$ have identical values.

\end{itemize}
\end{mydef} 

\begin{mydef} The process of finding  the quotient of a fracterm $P/Q$ is alternatively referred to as (i) calculating the fracterm $P/Q$, 
or (ii) as calculating the quotient of $P$ and $Q$, or (iii) as calculating the value of the fracterm $P/Q$.
\end{mydef}

There is considerable room for variation in these definitions. The notion of sameness of fracterms depends on a notion of 
sameness for numerators and for denominators. It is clear that difference of value implies non-sameness, but are the fracterms
$(2+3)/4$ and $(3+2)/4$ the same? Yes if one works modulo commutativity and no if one does not. 
Every notion of sameness for fracterms requires some additional, implicit or explicit, information about the particular use of that notion, 

\begin{mydef} 
\label{SimpleFr}
A flat fracterm is a fracterm of the form $P/Q$ where 
$P$ and $Q$ are closed or open terms over the signature of rings with naturals
in decimal notation.\footnote{%
In Dutch: ``vlakke breukvorm'' of ``eenvoudige breukvorm''. This understanding of simple fracterm deviates from common usage, 
as the traditional idea is that a simple fracterm is the same as a common fracterm.}
\end{mydef}

\begin{mydef} An inversion is an expression in the signature of inversive meadows or an extension thereof, of the form $P^{-1}$.
\end{mydef}
An inversion may be referred to as a fracterm with the understanding that $P^{-1} $ stands for $1/P$.
%

\begin{mydef} A flat form is an expression in the signature of meadows or an extension thereof which is not a fracterm and which has no
 fracterms or inversions as subterms.\footnote{%
In Dutch: ``vlakke vorm''.}
\end{mydef}

\begin{mydef} A mixed fracterm is an expression of the form $P + Q/S$ where $P$ is a flat form and $Q/R$ is a flat fracterm.
\end{mydef}
Using the notion of a flat form, a flat fracterm is a fracterm of which the 
numerator and the denominator are flat forms, and a mixed fracterm is a sum of a flat form and a simple fracterm. Flat forms involving addition and multiplication only are also referred to as polynomials.

\subsection{The fracterm quotient contrast in general}
The contrast between fracterm and quotient is present for other operations as well:
\begin{itemize}
\item $P+Q$ is the addition of $P$ and $Q$, the value of which is called the sum of $P$ and $Q$. 
$P+Q$ may be called a plusterm.
\item $P \cdot Q$ is the multiplication of $P$ and $Q$ the value of which is called the product of $P$ and $Q$. $P\cdot Q$ may be called a multerm.
\item $P-Q$ is the subtraction of $P$ and $Q$, the value of which is called the difference of $P$ and $Q$. $P - Q$ may be called a subterm.
\item $-P$ is the minus of $P$, the value of which is called the opposite of $P$. $-P$ may be called a minterm.
\end{itemize}
It is tempting to introduce a quotient operator as follows: $\quot{x}{y}$ denotes the quotient of $x$ and $y$. As a value the quotient has no structure, a quotient has no numerator and denominator and is not amenable to simplification. It follows that $x/y = \quot{x}{y}$ while $x/y \not\equiv  \quot{x}{y}$. However, doing so goes against the spirit of fracterms. By all means $P/Q$ is the way to denote the quotient
of $P$ and $Q$ and no purpose is served by having yet another notation for it.

\subsection{Secondary fracterm terminology}
\label{SecondaryFTT}
In the secondary fracterm terminology I provide an adaptation to fracterms of other aspects of the traditional terminology about fractions. For theoretical work most of the secondary terminology is immaterial, however. 
The following notions related to fracterms require an explanation: common fracterm (also vulgar fracterm), unit fracterm, proper fracterm, simplified fracterm, fracterm simplification, composite fracterm, and mixed fracterm. I refer to these notions as secondary fracterm terminology. From these notions simple fracterm and mixed fracterm are useful for theoretical work, and the definitions below work for open terms as well as for closed terms:

Now the other elements of secondary fracterm terminology are mainly used in school arithmetic and the some additional rules apply:
\begin{mydef} Assuming characteristic zero the following secondary fracterm terminology is reasonable/can be adopted:

\begin{enumerate}
\item 
A common fracterm (also called  vulgar fracterm)\footnote{%
I prefer to avoid the phrase simple fracterm (related to the well-known phrase simple fraction), 
because of the confusion that may arise, as with fractions, that simple fractions 
need not be simpilified.}
is a fracterm of the form $P/Q$ where $P$ and $Q$ are positive natural numbers in 
decimal notation without leading zeros.\footnote{%
In Dutch: ``gewone breukvorm''.}

\item 
An uncommon fracterm is a fracterm of the form $P/Q$ where $P$ and $Q$ are closed terms over the signature of rings with naturals
in decimal notation occurring a constants.\footnote{%
In Dutch: ``ongewone breukvorm''. It is not the case that a fracterm is either common or uncommon, for instance $(1/2)/3$ is neither.}

\item A closed problem fracterm is an uncommon fracterm $P/Q$ where the denominator satisfies $Q = 0$.
(The archetypical problem fracterm is $1/0$, and $(2+3)/(4-5) + 1)$ is a problem fracterm, $1/(x-x)$ is an open problem fracterm).
\item Two common fracterms $P/Q$ and $R/S$ are equivalent if $P \cdot S = Q \cdot R$. (It follows from this definition that
equivalent fracterms are equal; there is no other notation for fracterm equivalence than the equality sign, 
as in $P/Q = R/S$, which expresses about two fracterms just ``having the same value'', which is also the intended notion of equivalence).
\item Two common fracterms $P/Q$ and $R/S$ are the same if $P = R$ and $Q = S$.

\item A proper fracterm is a common fracterm for which the numerator is smaller than the denominator.

\item A unit fracterm is a proper fracterm for which the numerator equals $1$.

\item A fracterm is proper if it is common and the $\gcd$ of its numerator and its denominator is $1$.\footnote{%
In Dutch: ``echte breukvorm''.}
\item A fracterm is improper if it is common and not proper.\footnote{%
In Dutch: ``onechte breukvorm''.}
\item A {\em Scheinbruch} (German) is a common fraction of which the numerator is a multiple of the denominator.
\item A fracterm $P/Q$ is composite if either $P$ or $Q$ (or both) contains an occurrence of division.
\item A mixed common fracterm is an expression of the form $P~Q/R$ where $P$ is a natural number (expression) 
in decimal notation without any redundant leading zero and where $P/Q$ is a common fracterm.\footnote{%
In Dutch: ``gemengde echte breukvorm''.} (Every mixed common fracterm is also a mixed fracterm.)

\item A proper mixed fracterm is an expression of the form $P~Q/R$ where $P$ is a natural number (expression) 
in decimal notation without any redundant leading zero and where $P/Q$ is a proper fracterm.\footnote{%
In Dutch: ``gemengde breukvorm''.}

\item Two common fracterms $P/Q$ and $R/S$ are equivalent if $P \cdot Q = Q \cdot R$.
\item A fracterm $P/Q$ is simplified if it is proper and moveover $\gcd(P,Q)=1$.
\item A fracterm $P/Q$ allows simplification if it is proper and if it is not simplified.
\item A fracterm $P/Q$ is a simplification of  fracterm $R/S$ if both fracterms are common, both fracterms are equivalent, and
$P < R$.
\item Two fracterms $P/Q$ and $R/S$ have a common denominator (are like) if $Q$ and $S$ are the same. 
(Thus $2/(1+1)$ and $(1+1)/(1+1)$ are like but $(1+1)/2$ and $(1+1)/(1+1)$ are unlike (i.e. not like).)
\end{enumerate}
\end{mydef}
These definitions are somewhat flexible, for instance one may count $(-3)/4$ as a common fracterm, and one may count 
$-2~ 3/4$ as a mixed fracterm. Some authors avoid mixed fracterms altogether and, say,  insist writing $21+5/17$ rather than $21~5/17$

The relation between fracterm equivalence and fracterm equality is tricky. The idea is as follows: 
(i) that two fracterms which are not the same may be equivalent, 
(ii) fracterm equivalence can be defined without making use of the notion of the value of a fracterm, 
(iii) after some preparations it may be concluded that calling equivalent fracterms equal works in the same way as asserting that $2+2=1+3$. Both sides have the same value, but are not the same sums, 
(iv) for values the distinction between sameness and equality is not made.

I avoid the question whether or not the fracterms $(3+4)/5$ and $(4+3)/5$ are equivalent by 
having fracterm equivalence defined only for
simple fracterms (which these two are not).

The most well-known rule of calculation for fracterms is the so-called quasi-cardinality rule (QCR) which tels how to add two like
common fracterms thereby obtaining an uncommon fracterm with the same denominator. QCR is a phrase ascribed to Griesel~\cite{Griesel1981} in Padberg~\cite{Padberg2012}.

\begin{proposition} The sum $\alpha + \beta$ of two like common fracterms $\alpha \equiv P/R$ and $\beta \equiv Q/R$ 
 has the same value as the uncommon fracterm $(P+Q)/R$. In equational form: $P/R + Q/R = (P+Q)/R$.
\end{proposition}

\subsection{Arithmetical quantities, a container class for fracterms}
It is plausible to require an answer to the question: what kind of thing is a fracterm. Are fracterms a kind of some other type? 
The easy answer is that fracterms are terms, or expressions if one so prefers. 
Both replies have as a problem an implicit syntactic bias which I prefer to avoid. 
I prefer to use \emph{arithmetical quantity} (AQ) as the name of a container class name a container class of fracterms. The idea is that
$2+7$, $31 \cdot 56 -2 \cdot 3$ etc. are AQs, and if $P$ and $Q$ are AQs then so is $P/Q$. AQs of the form $P/Q$ are fracterms.

The introduction of AQ sheds light on some particular aspects of the design decisions which have been made thus far. Let $\Nat$ 
stand for the natural numbers and $\Int$ stand for the integers. The assumption that assimilation between names and named
objects takes place at the level of naturals and integers indicates that: 
$\Nat = |\Nat|$, $\Int = |\Int|$, and $|\Nat| \subseteq |\Int| \subseteq \AQ$.
Denoting the field of rational numbers with $\Rat$ and its domain with $|\Rat|$ in addition we expect $|\Int| \subseteq |\Rat|$.
Now rational numbers and fracterms have not been assimilated so that we may have $|\Int| = \AQ \cap |\Rat|$, 
with as a consequence that, say, $2/3 \notin |\Rat|$. I will say however, that $2/3 \in_= |\Rat|$, where $\_\in_=\-$ 
stands for: ``has the same value as an element of''. This $2/3$ is an AQ, it is a fracterm and it has the 
same value as some rational number (or is equivalent to some rational number). 
It is also the case that $2/3$ denotes a rational number. Now the question aries if it is also the case that $2/3$ is a number. 
I agree, in the context of fracterms as outlined above, 
with the assertion that ``$2/3$ is a number'' used as a shorthand for ``$2/3$ is a fracterm which denotes a number''. 
This shorthand works 
for a particular instance such as $2/3$ while it does not generalise to ``fracterms are numbers''.

\subsubsection{Selectors for numerator and denominator}
On AQ two relations exist with fracterms as domains: $Num$ and $Denom$ with the following definition:
\begin{mydef} $\mathrm{Num} \subseteq \AQ \times \AQ$, $\mathrm{Denom} \subseteq \AQ \times \AQ$ with the smallest graph so that for all fracterms $P/Q$:
$Num(P/Q,P)$ and $Denom(P/Q,Q)$.
\end{mydef}

The very notion of an AQ may be understood as an answer to the following question:" what kind of entity is supposed to come with a nominator and a denominator? The answer is  that  (i) some AQs but not all come with these two attributes, 
and (ii) no rational number comes with these attributes.

\subsection{Elements of a rationale for the definition of fracterms}
When choosing these descriptions of  terminology I had these requirements in mind:
\begin{enumerate}
\item Fracterms are not numbers because there is no plausible way in which numbers are equipped (decorated) 
with a numerator and a denominator
so that such decoration can be done in different ways for the same number.

Fracterms are also not (considered to be) ``a kind of numbers'' because it is not sufficiently clear what meaning to assign to that statement. 
\item In terms of assimilation and dis-assimilation (\cite{NicaudBG2001}): the default position is that the notions of 
fracterm and number are not assimilated. Clearly it is possible to work in a mode in which both notions are assimilated. I hold that
at any moment dis-assimilation of both notions ought to be an option, and that it is not the case that an assimilated perception of 
fracterms and numbers is more mature or more natural than a dis-assimilated perception of both notions. 
\item Speaking of addition of fracterms  suggests that there is something 
specific which can be said about the addition of fracterms in excess of the addition of AQs in general. 
In principle $P/Q+R/S$ is the sum of fracterms $P/Q$ and $R/S$. I will assume that
what is meant with ``addition of fracterms''  amounts to ``writing a sum of flat fracterms as a single flat fracterm'', 
and in the absence of variables ``writing a sum of common fracterms as a single common fracterm''. 
\item The idea of calculation cannot be explained or understood without having some notion of syntax at hand, 
for instance in the way it appears via the use of signatures: calculating an expression $P$ is finding in a systematic and defensible manner an expression 
$Q$ in a preferred form such that $P=Q$. Any notion of calculation has some parameters which are mostly left implicit: 
(i) a notion of equality (or equivalence), (ii) a preference for certain forms, and
(iii) conventions for writing instances of these preferred forms.
\item In order to define fracterms one need not  know (or define) in advance what rational numbers are.  
Moreover a particular definition of rational numbers, provides no substitute for a definition of fracterms.
\item Providing a systematic definition of fracterms is not a conventional way of going about the introduction of mathematical notions.
Fracterms are not a part of conventional mathematics. Like theorem, definition, proof and formula, fracterm is much used is but left undefined as a component of the informal language used for mathematical work. For notions which are somewhat 
syntactic mathematicians prefer not to provide definitions.
\end{enumerate}

\subsection{Problem fracterms}
A key notion is the idea of probem fracterms, which give rise to the issue of division by zero and ramifications thereof.
\subsubsection{Open (i.e. non-closed) problem fracterms}
Let $P/Q$ be an occurrence of a fracterm with variables among $x_1,\ldots,x_n$ in some text. 
In a mathematical exposition a specifica variable, 
say $x_i$ may occur as follows: (i) determined by a previous definition of the form ``let $x_i = t$'' with $x_i$ not in $t$, (ii)
conditionally existentially quantified, i.e. in the scope of a quantifier $\exists _{x_i : \phi} (\dots)$ 
(i.e. there is an $x_i$ such that $\phi(x_i)$ and $\ldots$), (iii) conditionally universally quantified
in the scope of a quantifier $\forall _{x_i : \phi} (\dots)$ (for all $x_i$ such that $\phi(x_i)$ it is the case that  $\ldots$.

A fracterm is a problem fracterm if, together with the information known about its variables, it is not possible, 
to the satisfaction of a reader, to  guarantee that $ \neq 0$.

\subsubsection{How to deal with problem fracterms?}
\label{HowToPF}
It is conventional to require of a text that it contains no closed or non-closed problem fracterms. Instead of this strict requirement, 
I propose that for an arithmetic it must be known how problem fracterms are to be treated.

\begin{description}
\item [Division by zero is impossible and must be prevented.] In a valid text there may not be any 
occurrence of a closed problem fracterm and there 
may not be any occurrence of a non-closed problem fracterm (prevention perspective). 

The reason for adopting a prevention perspective on division by zero lies in the observation that no
plausible extension of the rational numbers has ever been proposed in which $0 \cdot q = 1$ for some element $q$. The case 
differs from $\sqrt{-1}$ (which requires that $q^2+1 = 0$ for some $q$) which, although equally implausible at first sight, nevertheless has a plausible interpretation in various 
extensions of the rational numbers.

As a consequence of the prevention perspective a text may not even discuss the status of problem fracterms, at least not by means of 
discussing examples of such fracterms. This perspective is ``meadow based'' (and for that reason deviating from standard practice) as it acknowledges that, say $1/0$ is a fracterm, though
at the same time prevents it from being used. This perspective constitutes the best approximation of conventional practice within the 
setting of meadow based fracterm theory.
\item [Division by zero is legal and division is a partial function.] Closed problem fracterms are considered to have no value, 
and some chosen logic of partial functions  determine how truth values are determined of formulas which contain such fracterms.
This carries over to non-closed problem fracterms.

The logic of partial functions involved is a non-trivial parameter in this case. For instance ``$1/0 = 2/0$'' may be considered true, 
but may also be considered false or it may be considered to have no truth value 
(in which case some non-classical logic must be chosen). 
In each of these three cases two options are left for say $3/0 = 5$: false and no truth value. In all six cases ``$0/0 = 0/0$ is true''
is an option, and so is ``$0/0=0/0$ is false'' and ``$0/0=0/0$'' may have no truth value. If ``$0/0 = 0/0$'' is considered true, 
$0/0 = 0$ may or may not hold. The number of plausible logics for partial functions exceeds 10 at least.
\item [Division by zero is legal and produces zero.] All closed problem fracterms are set equal to $0$. 
For open fracterms substitution determines what happens. This is the preferred approach in the theory of meadows. 
\item [Division by zero produces unsigned infinity.] This is the preferred solution in the theory of wheels.
\item [Division by zero produces (positive) signed infinity.] The preferred solution in the theory of 
transrationals and transreals and other transfields (collectively known as transmathematics).
\end{description}

\subsection{Operations restricted to flat fracterms}
When working in the meadow of rationals $\Rat^d_0$, 
flat fracterms have pleasant closure properties w.r.t. the majority of arithmetical operations with the effect that several operations can simply be thought of as operations on fracterms:
\[-\frac{x}{y} = \frac{-x}{y}, \frac{x}{y} \cdot \frac{u}{v} = \frac{x \cdot u}{y \cdot v}, 
 \frac{x}{y} / \frac{u}{v} = \frac{x \cdot v}{y \cdot u}\] 
 For addition the situation is more complicated: addition is a function of type $\AQ \times \AQ \to \AQ$ which satisfies QCR
 \[ \frac{x}{z} + \frac{y}{z} = \frac{x+y}{z} \]
 However, addition does not restrict to a function of type $$\flatfracterm \times \flatfracterm \to \flatfracterm$$ 
 for instance there is no flat fraction $P/Q$
 such that $1/x + 1/y = P/Q$ (for a proof see~\cite{BergstraM2016}). 
 
 These properties together with $x \neq  \to x/x = 1$  imply CFAR, 
 the conditional fracterm addition rule:
 \[x \neq 0 \wedge y \neq 0 \to \frac{x}{y} + \frac{u}{v} = \frac{x \cdot v + u \cdot y}{y \cdot v}\]
 CFAR is a conditional which does not achieve the transformation of sums of flat fractions to flat fractions.
 
 In general, that is when working in other arithmetical datatypes than the meadow of rationals or the meadow of reals several of the 
 equations just mentioned may fail, and the status of flat fracterms becomes less central.
 
 If one restricts addition to common fracterms it is plausible to view addition as a relation. If one insists that addition restricts 
 to common fracterms as a function then using CFAR creates a discrepancy with QCR as different outcomes are obtained for $1/2 + 1/2$. The discrepancy with QCR can be settled by means of a definition based on  CFARcf, (the conditional fraction addition rule for 
 common fractions):
\[q \neq 0 \wedge s \neq 0 \to \frac{p}{q} + \frac{r}{s} \cong 
\frac{(p \cdot s + q \cdot r)\backslash \gcd(q,s)}{(q \cdot s)\backslash \gcd(q,s)}.\]
Here ``$\backslash$'' represents integer division and gcd produces the greatest common divisor.

\section{Fracterm collection oriented fraction definitions}
Having defined fracterms and fracterm related terminology as above it is possible to
specify different notions of a fraction. Besides providing at least one adequate definition of fractions (fractions as fracterms), 
surveying the plurality of fraction concepts constitutes an intended application of the introduction of fracterms.\footnote{%
The search for a systematic survey of fraction definitions constitutes an instance of vertical mathematising 
according to the classification proposed by Treffers in~\cite{Treffers1987}.}

Given a definition of fractions one may distinguish for each fraction (according to the definition) its base type and its target type.
The base type isha the fraction as an entity is, and the target type indicates what it refers to. base type and target type may be identical.

\subsubsection{Numerator Paradox}

The Numerator Paradox refers to a pattern of problematic arguments.
\paragraph{An instance of the Numerator Paradox:}
Consider the  fractions $2/3$ and $4/6$. The assertion $\mathrm{Num}(2/3,4) \wedge \neg \mathrm{Num}(2/3,4)$  will be derived
from some common assumptions. 

The (unique) numerator of $2/3$ is $2$, which is written as $\mathrm{Num}(2/3,2)$,
and the numerator of $2/3$ is not $4 $which is written as $\neg \mathrm{Num}(2/3,4)$. 
Now the numerator of $4/6$ is $4$ so that $\mathrm{Num}(4/6,4)$. For fractions $a/b$ and 
$c/d$ it is assumed that $a \cdot d = b \cdot c$ implies that both fractions are equal. It follows that $2/3 = 4/6$.
Now using $2/3 = 4/6$, substitution of equals for equals (2/3 for 4/6 in this case) 
in the assertion $\mathrm{Num}(4/6,4)$ yields $\mathrm{Num}(2/3,4)$. Taken together it follows that 
 $\mathrm{Num}(2/3,4) \wedge \neg \mathrm{Num}(2/3,4)$ as announced above. 
 
The latter assertion together with its  derivation constitutes an  instance of 
the Numerator Paradox. 

The Numerator Paradox consists of the entire family of assertions following the same pattern. 
 The Numerator Paradox can be formulated for each notion (definition) of fractions, 
and for each adequate definition 
of fractions it must be clear why the Numerator Paradox does not apply.
\begin{proposition} Fracterms are not numbers.
\end{proposition}
\begin{proof} If fracterms are numbers the Numerator Paradox applies, which must be rejected.
\end{proof}

\subsection{Five single fracterm based definitions of a fraction} 
A definition of fractions in terms of fracterms is called single fracterm based if according to the definition 
each fraction has a unique fracterm as its counterpart which represents the mathematical entity which stands for that fraction. 
\subsubsection{Fractions as (frac)terms}
As a first definition of fractions I mention the definition which I prefer, which constitutes the fractions as terms paradigm on fractions.
\begin{mydef}
\label{Fraction}
A fraction is a fracterm.
\end{mydef}
In this case the Numerator Paradox is easily resolved: 
it is wrong to assume that $\mathrm{Num}(P/Q,T) \wedge P/Q = R/S$ implies 
$\mathrm{Num}(R/S,T) $. What holds instead is: $\mathrm{Num}(P/Q,T) \wedge P/Q \equiv R/S$ implies 
$\mathrm{Num}(R/S,T) $.

Upon adopting this definition of fractions, the fraction base type is fracterm, and the fraction target type is rational number.
If one adopts this definition, and if one is willing to do without other perspectives on the notion of a fraction,
the word ``fracterm'' may be eliminated from the discussion altogether and in
the list of definitions above everywhere fracterm may be replaced by fraction.

The list of definitions related to fracterms 
has been prepared in such a manner that the list resulting from the replacement of fracterm by fraction produces 
a fully workable fraction terminology.\footnote{%
Corresponding with replacing fracterm by fraction in the Dutch translations, as given in various footnotes ``breukvorm''
 is to be replaced by ``breuk''.}
 \paragraph{Doubts about fractions as fracterms.}
I advocate fractions as fracterms as the most workable definition of fractions for theoretical work. 
It has been used in~\cite{BergstraT2008,BBP2013, BergstraM2016}, and its fits easily in a context where logic and data type theory are prominent. But disadvantages of defining fractions as fracterms must be taken into account as these disadvantages serve as the 
core of motivation for looking for other definitions. 

First of all one may disagree for various reasons with the idea that fracterms as a concept can be integrated in smooth explanation of elementary arithmetic. Then, assuming that one is satisfied with the definition and use of fracterms, still the following objections against taking fractions for fracterms can be made.
\begin{enumerate}
\item the notion of a fraction remains fraught with lack of clarity. For instance  is $(1+1+1)/1$ a fracterm? If so why is it possible to avoid brackets in $1+1+1$ without having an extended description of working modulo associativity and commutativity, if not, is that workable? 
\item The notion of being the same is not sufficiently clear for fractearms: do we expect $1/(x+(y+z)) \equiv 1/((x+y)+z)$? 
If so, on what basis, if not should the concept of fraction be made more abstract? 
 Is $1+1$ a fraction, if so why, if not, what distinction to make between $1+1$ and $(1+1)/1$?
\item Given that $2/3$ is the most concise way to denote ``the number denoted by $2/3$'', why not say that $2/3$ is a number. This 
constraint is counterintuitive.
\item There is remarkably little support for fractions as fracterms in the literature on fractions, which matters for whoever deal with
 fractions in a science based manner.
\end{enumerate}

\subsubsection{Fractions as division safe fracterms}
Following~\cite{BergstraT2008} a fracterm is division safe if it is not a problem fracterm and none of its subterms are problem fracterms. 
This notion works for open expressions under the assumption that concerning the variables of a fracterm 
enough information is available to
infer that its denominator is non-zero.
\begin{mydef}
\label{Fraction}
A fraction is a division safe fracterm.
\end{mydef}
Upon adopting this definition of fractions, the fraction base type is fracterm, and the fraction target type is rational number.
Upon adopting this concept of fractions the list of definitions may be adapted by replacing  fracterm by fraction throughout and 
adding the requirement that a denominator is non-zero sufficiently often to be able to make that inference whenever a fracterm is used.

The notion of a problem fraction will not arise as problem fracterms don't qualify as fractions under this convention. 
All other notions regarding fracterms can be carried over to fractions, however.

\subsubsection{Fractions as values}
Very common is the perspective that  a fraction is a rational number. 
In fact this perspective on fractions is independent of one's conception of rational numbers. 
However, in the context of the current survey I will assume that rational numbers are defined as in Paragraph~\ref{Rationals}.
\begin{mydef}
\label{FractionAsValue}
A fraction is a value, in particular a fraction is a rational number, which corresponds in the setting of the current paper  
to: a fraction is a simplified common fracterm. 
\end{mydef}
Upon adopting this definition of fractions, the fraction base type is rational number, and the fraction target type is rational number.
In the fractions as values paradigm the Numerator Paradox does not apply because the relations $\mathrm{Num}$ and $\mathrm{Denum}$ don't exist on values. 
\subsubsection{Fractions as fracpairs I: pairs of integers}
\begin{mydef}
\label{Fraction}
A fraction is pair $(n,m)$ with $n,m \in \Int$ and $m \neq 0$. In the current setting a fraction is a common fracterm.
\end{mydef}
Upon adopting this definition of fractions, the fraction base type is fracpair, and the fraction target type is rational number. The Numerator Paradox is resolved by noticing that substituting equals for equals does not work in the predicates Num and Denom.
\subsubsection{Fractions as fracpairs II: pairs of rationals}
\begin{mydef}
\label{Fraction}
Upon adopting this definition of fractions, the fraction base type is pair of rational numbers, and the fraction target type is rational number.
A fraction is pair $(p,q)$ with $p,q \in \Rat$ and $q \neq 0$. In the current setting a fraction is a  fracterm for which both the numerator and the denominator are simplified common fracterms, and for which the denominator is not equal to zero.
\end{mydef}

\subsection{Differential fraction talk}
The five definitions above lead to different properties of fractions: for instance if fractions are considered division safe fracterms, values, or pairs the notion of a problem fraction evaporates. Indeed, the notion of a problem fracterm, and the corresponding ramification of 
ways for dealing with problem fracterms only appears in the setting of fractions as fracterms. All other primary and secondary terminology
works for fractions as division safe fracterms too.

If fractions are considered values, pairs of integers or pairs of rationals the notion of a composite fracterm does not apply. If fractions are considered values then the notion of simplification of fractions disappears. If fractions are considered pairs of rationals then some fractions are simplified but not all fractions admit simplification and the the notion of fraction simplification does not apply.

\subsection{Beyond the five single fracterm based fraction definitions}
Besides of the five fracterm oriented fractions, however, several other options for defining fractions can be found. I provide an incomplete survey of six additional options for defining a fraction in the following Paragraphs, followed by a parametrised family of $26$ fraction definitions in Paragraph~\ref{Ufracterms}.
\subsubsection{Numbers as sets of fractions}
This survey is incomplete in the sense that certainly not all conceivable views on fractions are covered. For instance on may insist that 
rational numbers are different objects than what has been assumed above. 
One may insist that rather than a pair $(n,m)$ of integers with $m$ positive and $n,m$ relative prime, a rational number is a collection
of pairs $(n,m)$ with $m \neq 0$ which are pairwise equivalent (with equivalence given by $(n,m) \equiv (a,b) \iff n \cdot b = a \cdot m$).
Now one may insist that the following assertion is true for one's preferred notion of fraction:
\begin{quote} 
 $\Psi_{cf/rn} $: A common fraction is an element of a rational number.
\end{quote}
Holding  $\Psi_{cf/rn} $ is compatible with assimilating say $(2,3)$ and $2/3$ while not assimilating $(2,3)$ with $4/6$.
The statement $\Psi_{cf/rn} $ is not valid for any of the perspectives that are covered in this paper. If numbers are understood as sets 
(i.e. equivalence classes) of common fracterms the following candidate definition of fractions arises:

\begin{mydef} A fraction is an element of a natural number.
\end{mydef}
\subsubsection{Operator view}
The operator view of fractions as used in~\cite{AthenG1978} is an example of another approach not covered in this paper. in the operator view a fraction is understood as an operator ($\lambda x. P \cdot x$) which can gradually become known by inspecting examples for various values of $x$. 

\begin{mydef} A fraction is a linear mapping from rational numbers to rational numbers which takes value zero on zero.
\end{mydef}
The Numerator Paradox fails on the fact the Num and Denom don't exist as relations of linear functions from rationals to rationals.
\subsubsection{Fracvalues, fracpairs, and fracsigns}
Besides fracterm there is room for other new terminology: fracpair, fracvalue, and fracsign. A fracvalue is the value of a fracterm, 
a fracsign is a sign which refers to a fracterm, and a fracpair is pair which contains the 
values of the numerator and the denominator of a fracterm.

In this paper fracvalues are the rational numbers, i.e. the elements of $|\Rat|$. Fracsigns are an open ended class of physical entities.
When introducing fracsigns different fracsigns may refer to the same fracterm. Fracpairs are sensitive to the type of the fracterm components. If one is dealing with non-composite fracterms only it is plausible that fracpairs are chosen as pairs of integers, 
perhaps with restrictions on the second member of the pair (nonzero, non-negative, exceeding 1 etc.).
If one contemplates composite fractions the restriction that fracpairs are pairs of integers is to strong, and pairs of rational numbers may be needed.

Returning to the matter of incompleteness: a perspective on fractions where fractions are understood as fracsigns of some kind 
falls outside the class of interpretations of fractions as a kind of fracterm. 

\begin{mydef} A fraction is a fracsign.
\end{mydef}
With this definition $1/2$ and $1/2$ (different occurrences of fracsigns for  the same fracterm) are not the same. This definition also calls for natural numbers and integers to be identified with appropriate signs. In any case it is implausible that signs are identified with numbers.
\begin{mydef} A fraction is the referent of a fracsign.
\end{mydef}
\subsubsection{A multi-fracterm definiton of fractions I}
Fractions may be different kinds of fracterms. For instance:
\begin{mydef} A fraction is either a pair $(\mathsf{fp},P)$ with $P$ a common fracterm (i.e. a fracpair) or a pair $(\mathsf{r},P)$ with $P$ a simplified common fracterm (i.e. a rational number).
\end{mydef}
In this definition both the fracterm $(4,6)$ and the number  $(2,3)$ are simultaneously assigned to $4/6$.
Some fractions are values and don't have a numerator and a denominator, while other fractions are decomposable into a 
nominator and a denominator and are not values. This definition of fractions supports polymorphism. The point of departure
is fracsigns $P/Q$ which may refer, depending on the context to either one of two fracterms.
 It depends on the context how a fracsign $P/Q$ is understood. When one reads/writes, ``the fraction $9/12$ has an odd numerator'', 
 the fracsign  $9/12$ is meant to be  the fracpair $(9,12)$ (or equivalently  the common fracterm $9/12$), 
 while when one reads/writes $9/12>1/2$ the fracsign $9/12$  the 
 refers to the fracvalue $(3,4)$ (or equivalently the simplified common fracterm $3/4$). In this case composite fractions do not exist (i.e. such eniities are not fractions).
\subsubsection{A multi-fracterm defintion of fractions II}
Another multi-fracterm definition of fractions is:
\begin{mydef} A fraction is either a pair $(\mathsf{ft},P)$ with $P$ a  fracterm or a pair  $(\mathsf{r},P)$ with $P$ a simplified common fracterm (i.e. a rational number).
\end{mydef}
In this context it is valid to say that ``the composite fraction $1/(2/3)$ equals the number $3/2$'' and also that ``the number $(1/2)/3$ equals a unit fraction.'' Moreover ``some fractions are rational numbers, and all rational numbers are fractions'' (with intended meaning: 
``some fracsigns refer to rational numbers and  all rational numbers are referred to by a fracsign in an appropriate context'').

\subsubsection{Indeterminate fraction definitions}
\label{Ufracterms}
Let $U$ be a subset of $\{\mathsf{ft},\mathsf{dsft},\mathsf{fv},\mathsf{fpi},\mathsf{fpq}\}$, and assume that $U$ has at least two elements
(there are $26$ such options for $U$). 
These codes indicate the respective five fracterm oriented fracterm definitions. Now ``fractions as $U$-fracterms'' works as follows:
\begin{mydef} A fraction is understood as a fraction in the sense of $\alpha \in U$, for some $\alpha$ with the understanding that for each context which is dealing with fractions a unique $\alpha$ must be chosen.
\end{mydef}
The notion of a context is not entirely precise but a new context is started in a text if no previously introduced names play a role anymore.
Typically that is the case if in a textbook a new Paragraph starts. 

A remarkable property of fractions as $U$-fracterms is that the logic becomes somehow non-deterministic while it need not be possible, as in the polymorphic case, to infer from the context which of the specific interpretations of fractions is used. That must perhaps be discovered by inspecting the use. As a consequence it may be the case that the logic of fractions resembles an informal logic rather than  a formal logic. This consequence of underspecification of the notion of a  fraction may be reinforced if additional degrees of freedom are
introduced by admitting different policies for assimilation and dis-assimilation in connection with naturals and integers. These matters leave room for further research.

\subsection{From fracpairs to fracvalues by way of chunk and permeate}
In Bergstra \& Bethke ~\cite{BergstraB2015} a fractions as pairs perspective to the logic of fractions 
is followed using an application of the  chunk and permeate method of 
Brown \& Priest~\cite{BrownP2004} which, like the Logic of Paradox LP, 
is considered to belong to the realm of paraconsistent reasoning. 
 From the theory of positive rationals all the equations as well as
 the inequation $1 \neq 2$ are permeated to the new theory, thereby forgetting (chunk) 
 all information about numerator and denominator, finally the equation $x/x = 1$ is 
 added and then a workable theory of positive rational numbers emerges.

\section{Literature on fractions}
The literature on fractions contains at least 5000 papers. I have been unable to develop a systematic survey of views on fractions as 
put forward in these works. This section provides a brief discussion of work on what is a fraction and of papers which contain claims or assertions which I understand as being close to one of the three paradigms: fractions as terms, fractions as values, and fractions as pairs.
\subsection{Literature on the question ``what is a fraction?''}
 Van Hiele in~\cite{Hiele1973} suggests to do away with fractions entirely and only to use the inverse function instead. 
 Van Hiele considers fractions harmful.
 Schippers \cite{Schippers2014} provides a comprehensive survey of the history and development of fractions. 
 
Olanoff, Lo \& Tobias~\cite{OlanoffLT2014} use the phrase fraction content knowledge which suggest that 
fractions transcend a precise definition.
 
Padberg~\cite{Padberg2012}  assumes that fraction is a complex notion amenable to a thematic decomposition. Decomposing the notion of a fractions in so-called subconstructs originates from Kieren~\cite{Kieren1976}. Different subconstructs go with a different langue and notation and even a different way of thinking, thus constituting  different logics of fractions so to speak. 
Tanton~\cite{Tanton2014} (page 3) writes concerning the question   
``What is a fraction?'':
\begin{quote}
As we shall see answering this question is far from easy! It is incredibly hard to pin down exactly what a fraction is.
\end{quote}
Quinn~\cite{Quinn2013} indicates that fractions transcend the particular values they 
denote in a specific mathematical structure. Fractions are said to be procedurally generic. This observation is
convincingly linked to modern informatics. The quote explains an aspect which, however, I have not been able to 
take into account.
\begin{quote}
The net result is that for elementary use it is unnecessary to say anything about what or where a fraction is, 
as long as we avoid dividing by zero.\\

The corresponding drawback is that attempting to explain the ``what'' and ``where'' of fractions is likely to turn a straightforward procedural subject into a confusing mystery. The precise details show that writing a fraction specifies an object in the image of every appropriate ring morphism, not just one. One might think of this as being like a picture on the internet: it doesn't appear everywhere; one must go to a device with an internet connection, but it appears in every such context. This is a perfectly functional way to think about internet pictures. A full-precision description as functions, from internet-connected devices and URLs to images, would not help anybody. Falsehoods analogous to the ancient descriptions of fractions are not helpful either because they interfere with development of genuinely functional understanding.
\end{quote}
The following quote, taken from the Paying Attention to Mathematics series\footnote{%
\url{http://www.edu.gov.on.ca/eng/literacynumeracy/LNSAttentionFractions.pdf} (accessed February 26, 2017.}~highlights a popular style of providing ``definitions'' which I find disappointing:
\begin{quote}
What Is a Fraction?\\
\newline
A fraction most simply represents a number, as shown on the number line below:\\
(pictorial display of number line)\\
Yet within this very simple description lie some highly complex mathematical constructs that will be explored in this document. These constructs --- part--whole relationships, part--part relationships, quotient and operator --- are not mutually exclusive; they are different ways to represent and think about fractions.
\end{quote}  
This answer does not state what a fraction is but rather what it does, i.e. what it is meant to be used for. 
In~\cite{Strother2016} Strother et. al. propose 5 key ideas for teaching fractions. One of the key ideas is presented as follows:
\begin{quote}
\emph {Use Precise Definitions of the Numerator and Denominator}\\
A common definition in U.S. textbooks for the numerator and denominator is that the denominator is the parts to make a whole and the numerator is the number of parts being counted (Lamon 2012). However, these definitions do not adequately represent the way fractions must be understood to extend students' learning beyond the simplest examples of these numbers nor are these definitions sufficient to apply to learning fraction and decimal operations (Wu 2011). Reviewing the history of these concepts reveals that the denominator is intended to identify, or name the unit being counted and therefore the size of this unit.$\ldots$
\end{quote}

Clarke, Roche \& Mitchell in~\cite{ClarkeRM2011} emphasise the need for a holistic approach to fractions:
\begin{quote}
The dilemma for both teachers and students is how to make all the appropriate connections so that a mature, holistic and flexible understanding of fractions and the wider domain of rational numbers can be obtained, because fractions form an important part of middle years mathematics curriculum, underpinning the development of proportional reasoning, and are important for later topics including algebra and probability.
\end{quote} 
Redmond~\cite{Redmond2009} writes:
 \begin{quote}
 The research above shows the importance of developing fraction sense before introducing the rules and procedures for fractions operations and described ways to help
students make sense of fractions. However, the teaching and learning of division of fractions is an even more complex issue.
 \end{quote} 
 The notion of fraction sense is not defined in a systematic manner, but I assume that developing fraction sense 
 is subsumed under the mentioned holistic approach. 

The phrase ``paradoxical character of mathematical knowledge'' appears in Duval~\cite{RDuval2000}. Although Duval's paper is not dealing with 
fractions per se, his work is to the best of my knowledge unique within the educational literature on mathematics in that it explicitly highlights the paradoxical character of 
mathematical knowledge as a consequence of the presence of different representations of mathematical entities, 
the key phenomenon governing the relation between fractions and numbers, an issue already appearing with the seemingly trivial gap between binary representation and decimal representation of natural numbers. In his words:
\begin{quote}
Concerning the cognitive mode of access to objects, there is an important gap between mathematical 
knowledge and knowledge in other sciences such as astronomy, physics, biology, or botany. 
We do not have any perceptive or instrumental access to mathematical objects, 
even the most elementary, as for any object or phenomenon of the external world. 
We cannot see them, study them through a microscope or take a picture of them. The only way of gaining access to them is using signs, words or symbols, expressions or drawings. But, at the same time, mathematical objects must not be confused with the used semiotic representations. This conflicting requirement makes the specific core of mathematical knowledge. And it begins early with numbers which do not have to be identified with digits and the used numeral systems (binary, decimal).
\end{quote}
Moreover and most important for the case of fractions Duval highlights the importance of conversion between different representation, stating that ``conversion depends on incongruence''. 

 Van den Heuvel-Panhuizen~\cite{VandenHeuvel2003} investigates the development of bar models for fractions. 
 Bar models belong to the fraction construct ``part of a whole''. The educational value of such models is claimed to reside not so much in the models per se, and not in the objective of finding a best and definitive model as a tool for the student. Rather:
\begin{quote}
Formulated more precisely, it is not the models in themselves that make the growth in mathematical understanding possible, but the students' modelling activities. 
\end{quote}

In Lortie--Forgues, Tian \& Siegler~\cite{LortieTS2015} it is noticed that:
\begin{quote}
A fraction has three parts, a numerator, a denominator, and a line separating the two numbers. 
This configuration makes fraction notation somewhat difficult to understand. For instance, students,
especially in the early stages of learning, often misread fractions 
as two distinct whole numbers (e.g., 1/2 as 1 and 2), as a familiar arithmetic 
operation (e.g., 1 + 2) or as a single number (e.g., 12) ......
\end{quote}
According to these authors not only is a fraction not a single number, it is not (like the result of) a familiar operation on numbers either, and it is not  a pair. These remarks provide ample evidence of the difficulty of defining a fraction
in unambiguous positive terms.

Bruin--Muurling~\cite{Bruin2010} (p17) states:
\begin{quote}
What makes the concept of fraction so
difficult is that $\frac{a}{b}$ is both the division $a \div b$ and the factor between these two numbers, i.e. 
it expresses a relationship between two quantities and it defines a new quantity.
\end{quote}
This formulation definitely expresses that a rational number is missing some aspect(s) of a fraction. It also
acknowledges that the concept of a fraction is difficult.

\subsection{Support from the literature for fractions as (frac)terms}
Although viewing fractions as terms is a very obvious perspective on fractions for someone with a background in logic or computer science, 
I have not been able to find much support for it in the literature on fractions. 
 
 The claim that the fractions as fracterms perspective on fractions is meadow based is not entirely 
 precise and labeling this perspective as 
 meadow signature based is more adequate. 
The word ``fracterm'' is not used in the literature. 
I wil include under the fractions as fracterms subconstruct each approach where a fraction is implicitly or explicitly identified with an expression. I must acknowledge that some quotes (and my interpretation thereof) 
may in fact not represent the view of the authors as expressed in their papers at large. 

A fairly explicit proposition that fractions are to be viewed as terms, and that operations on 
fractions for that reason are operations on terms as well as on numbers, can be found in Mueller~\cite{Mueller1961}. 
Mueller uses symbol for term and does not explicitly distinguish between symbols and signs. Mueller's paper was written in reaction to 
a  proposal by Van Engen in~\cite{Engen1960} who took a fractions as fracpairs position.

Your Dictionary\footnote{%
\url{http://www.yourdictionary.com/fraction} (accessed February 25, 2017).} states that
\begin{quote}
The definition of a fraction is a mathematical expression with a numerator and a denominator, a disconnected piece or a small part of something.
\end{quote}
On ``What is a fraction''\footnote{%
See \url{http://www.cut-the-knot.org/WhatIs/WhatIsFraction.shtml} (accessed Februari 25, 2017).} 
it is stated that:\footnote{%
On Encyclopedia of Mathematics, however one reads (\url{https://www.encyclopediaofmath.org/index.php/Rational_number} accessed February 27, 2017) that a rational number (is) 
A number expressible as a fraction of two integers. Combining both definitions leads to a circularity which would be absent
in the presence of an independent notion of expression on which fractions may be based.}
\begin{quote} 
A \emph{fraction} (sometimes a \emph{common fraction}) is a way of expressing a number that is a ratio of two integers.
\end{quote} and also:
\begin{quote}
The number expressed by a fraction has infinitely many fractional representations.
\end{quote}
In this context I assume that a way of expression is the same as an expression. Fractions as terms may be hard to distinguish from fractions as names. For instance on 
Math Open Reference\footnote{%
See \url{http://www.mathopenref.com/fraction.html} (accessed Februari 28, 2017)} the following definition appears:\begin{quote}
A fraction is two quantities written one above the other, that shows how much of a a whole thing we have. For example we may have three quarters of a pizza:

\quad $\displaystyle \frac{3}{4}$ of a pizza 
The description reads as if an expression is specified but lacking more detail concerning the components of the expression this fragment may also be considered representative for fractions as names. 
\end{quote}
Wolfram MathWorld\footnote{%
See \url{http://mathworld.wolfram.com/Fraction.html} (accessed February 26, 2017). 
The original uses colours instead of italics for highlighting purposes.}
defines a fraction as follows:
\begin{quote}
A \emph{rational number} expressed in the form $a/b$ (in--line notation) or $\displaystyle \frac{a}{b}$ (traditional "display" notation), where  $a$ is called the \emph{numerator} and $b$ is called the \emph{denominator}. 
\end{quote}
I assume that a rational number expressed in a certain form is equivalent to: an expression of a certain form
which is supposed to denote a rational number. This definition provides a hybrid between 
fractions as expressions and fractions as numbers. In a fractions as expressions interpretation one misses the option to understand the fraction as denoting an element of a finite field. In a fractions as values interpretation of this definition the  expression comes as and additional attribute which is commonly not implicit in the 
notion of a rational number.

I consider fracterms to constitute a form of names. Viewing fractions as names imposes fewer constraints than viewing fractions as fracterms. 
Fractions as names  is central to the work of Rollnik in~\cite{Rollnik2009}.

The Merriam Webster dictionary provides several definitions of fraction, the one relating to arithmetic reads as follows:
\begin{quote} a numerical representation (as $\sfrac{3}{4}, \sfrac{5}{8}$, or $3.234$) indicating the quotient of two numbers.
\end{quote}
Assuming that $R$ being a representation indicating $s$ is comparable to $R$ naming $s$ it is fair to classify this definition as belong to the fractions as names paradigm. On Math Goodies Glossary\footnote{%
See \url{http://www.mathgoodies.com/glossary/term.asp?term=equivalent\%20fractions} (accessed February 25, 2017.)} equivalent fractions are thus specified:
\begin{quote}
Equivalent fractions are different fractions that represent the same number.
\end{quote}
In Math.com\footnote{%
See~\url{http://www.math.com/school/subject1/lessons/S1U4L1DP.html} (accessed on Februari 20 2017.)} it is asserted that:
\begin{quote}
There are many ways to write a fraction of a whole. Fractions that represent the same number are called equivalent fractions. This is basically the same thing as equal ratios. For example, $\sfrac{1}{2}$, $2/4$, and $4/8$ are all equivalent fractions. To find out if two fractions are equivalent, use a calculator and divide. If the answer is the same, then they are equivalent.
\end{quote}

Here rational numbers (also referred to as ratios) 
are identified with decimal expansions to which fractions can be evaluated, for instance by means of a calculator. Fractions are considered representations of numbers.
On TeacherVision\footnote{%
See \url{https://www.teachervision.com/operations/finding-equivalent-fractions-simplest-form} (accessed February 26, 2017.}
fractions are said to  name an amount:
\begin{quote}
Equivalent Fractions: Fractions that reduce to the same number and have an equal value. They are fractions that name the same amount in different ways.
\end{quote} Counting amounts as values this corresponds to a fractions as names position.
Wong \& Evans in~\cite{WongE2011} state
\begin{quote}
Understanding fraction equivalence necessitates students recognise that two or more fractions can represent the same quantity, thus belonging to an equivalence set.
\end{quote} I consider texts where fractions are considered representations of numbers or representations of quantities to belong to the same logical category as texts speaking of names.

Erdem~\cite{Erdem2016} provides a survey of fractions an insists that say $2/3$ has two meanings, 
a rational number and a fraction which must not be confused. Erdem also insists that on conceptual grounds
a fraction cannot be negative so that $-2/3$ can only be a rational number. Erdem states that 
\begin{quote}
In the broadest sense of the term, fractions are the expressions that are used to represent rational numbers and that can exist in infinite numbers.
\end{quote}
 I consider his view as a fractions as names view, rather than as fractions as terms view, because no mention is made of a specific syntax for 
fractions. He cites Vamvakoussi \& Vosniadou~\cite{VamvakoussiV2010} as follows: 
\begin{quote}
One could say that $0.5, 500/1000$ and also $1/2$ or $7/14$, are all rational numbers. However, a closer examination shows that they all have the same value, thus the accurate thing to say would be that these are alternative representations of the same rational number. A mathematician would define this rational number as the equivalence class of pairs $[a,b]$ such that $2a = b$, and would say that all the above are representatives of this class.
\end{quote}

Pot (\cite{Pot2005}) lists several problems and ambiguities with the term fraction and he proposes to use the term ``breuktal'' (pseudo Dutch, translating to fractional number) instead and to view that as a rational number. Given 
$r = p/q$ the numerator and denominator are taken to be to be the factors of $p$, resp. $q$ that are obtained after simplification. This approach qualifies as a fractions as values paradigm. Pot introduces an additional phrase (breuktal), thus preserving for fraction (breuk)  the option to denote something else than number.

\subsection{Support from the literature for fractions as values}
\label{FaVlit}
Rollnik in~\cite{Rollnik2009} provides a proposal for thinking about fractions based on the view that a fraction is a number, in particular a rational number. A thorough argument for the necessity of being able to give (and teach) a proper definition of fractions is given by Greenleaf in~\cite{Greenleaf2006}, who follows the  fractions as values paradigm.

Filep~\cite{Filep2001} writes:
\begin{quote}
The next step in the development of the fraction concept was its comprehension as a relative quantity to the whole as unit. The final step was made by Omar Khajjiam (XI. century) who interpreted it as fractional number arising from the division of two numbers.
\end{quote}
The Oxford Dictionary\footnote{%
\url{https://en.oxforddictionaries.com/definition/fraction} (accessed February 25 2017.)} provides as its first definition of a fraction:
\begin{quote}
A numerical quantity that is not a whole number (e.g. $1/2, 0.5$).
\end{quote}
However, the Oxford dictionary online doesn't claim that fractions have a numerator and a denominator. 
Instead those are said to be components of a vulgar fraction, and a vulgar fraction is a fraction expressed by a numerator and denominator (that is, not as a decimal fraction).
In Dictionary.com, however a fraction is also considered a number while a numerator is understood as follows:\footnote{%
\url{http://www.dictionary.com/browse/numerator?s=t} accessed, March 7, 2017.}
\begin{quote}
\emph{Arithmetic.} the term of a fraction, usually above the line, that indicates the number of equal parts that are to be added together; the dividend placed over a divisor:\\
\emph{The numerator of the fraction 2/3 is 2.}
\end{quote}
Usiskin uses in~\cite{Usiskin2007}  the phrase ``simplification of rational numbers'' instead of the more common 
phrase ``simplification of fractions''. From this use of terminology I infer that Usiskin views fractions as rational numbers.

\noindent In Bruce~\cite{Bruce2014} it is stated under the heading of ``multiple meanings of fractions'' that
\begin{quote}
A fraction is one quantity or amount
\end{quote}
Precisely this view corresponds to fractions as values which constitutes the core of this paper. 
A similar position is found in Fazio \& Sigler~\cite{FazioS2011}:
\begin{quote}
Fractions are numbers
\begin{quote}
Students need to understand that fractions are numbers with magnitudes.
\end{quote}
Research findings\\
\newline
Fractions are often taught using the idea that fractions represent a part of a whole. For example, one-fourth is one part of a whole that has been split into four parts. This interpretation is important, but it fails to convey the vital information that fractions are numbers with magnitudes....
\end{quote}
In ``Fractions as Numbers''~\cite{NationalC2014} it is asserted under the 
heading ``Conceptual Understanding''
 that:
\begin{quote}
Develop understanding of fractions as numbers, such as the following:\\
\newline
1/5 is the same as 2/10.\\ 2/5 is the same as 4/10.\\ 3/5 is the same as 6/10.
\end{quote}

\noindent Tucker~\cite{Tucker2008} states that
\begin{quote}
As soon as is judged feasible, a learner should be given the following definition of a fraction. 
This is the definition of a fraction used in many other countries:\\
\newline
A fraction is a number that is an integer multiple of some unit fraction.\\
\newline
In mathematical notation, we mean a number of the form $k(1/l)$, for whole numbers $k, l (l > 0)$. This definition assumes that the person has first developed a good understanding of what a unit fraction is. Unit fractions are discussed extensively in the next section. Note that in the essay, we will not worry about more complicated fractions, with numerators and denominators that are themselves fractions or irrational numbers. 
\end{quote}
In Wu~\cite{Wu2014} it is stated on 
page 19 that
\begin{quote}
The fundamental fact states that a fraction is not changed (i.e., its position on
the number line is not changed) when its numerator and its denominator are both
multiplied by the same nonzero whole number.
\end{quote}
The notions of fraction equality and fraction equivalence are identified. From the quoted statement 
and several other assertions in~\cite{Wu2014} I conclude that fractions as values is adopted as the 
underlying view of  fractions in this work. Tanton~\cite{Tanton2014} writes
\begin{quote}
So here is the best answer I have to the question: what is a fraction?
\begin{quote}Fractions are some kind of numbers that satisfy the basic five beliefs stated above (and satisfying all the usual rules of arithmetic too). Everything you want to do and understand about fractions follows from these five rules.
\end{quote}
Let's see now how everything does indeed unfold from these five belie(f)s.
\end{quote}
The mentioned beliefs are common rules of calculation on positive 
rational numbers expressed with addition, multiplication, and division.

A remarkably systematic site on fractions is provided by Skillswise.~\footnote{%
\url{http://ashleyroad.aberdeen.sch.uk/wp/wp-content/uploads/2014/04/FFractions.pdf}, accessed on April 7 2017.}
Only simple fractions are discussed and equivalent fractions are called the same, 
with the consequence that different fractions can be the same. I conclude that the site takes fracvalue as the fraction base type. However, a different interpretation is possible: the same fractions $R$ and $S$ can be different unless 
$R$ and $S$ are identical. If sameness is understood thus, the fraction base type used in this site is best identified as 
 fracpair rather than as fracvalue.
 
In a conceptual paper on the difficulty of elementary mathematics Wu~\cite{Wu2009} states (p8):
\begin{quote}
... understanding fractions is the most critical step in understanding rational numbers and in preparing for algebra. 
\emph{In order to learn fractions, students need to know what a fraction is.} Typically, our present math education lets them down at this critical juncture. All too often, instead of providing guidance for students? first steps in the realm of abstraction, we try in every conceivable way to ignore this need and pretend that there is no abstraction.
\end{quote}
Wu then proceeds to define fractions as rational numbers which are considered point so  the number line, 
thereby leaving unexplained what 
numerators and denominators are. The following assertion (p8 in~\cite{Wu2009}), however, casts some doubts on the status of fractions: are these rational numbers or not?
\begin{quote}
In particular, regardless of whether a number is a whole number, a fraction, a rational number, or an irrational number, it takes up its natural place on this line. 
\end{quote}
More specifically Wu states (p12):
\begin{quote}
In a formal mathematical setting, we now use this particular point as the official representative of $\sfrac{1}{3}$. In other words, whatever mathematical statement we wish to make about the fraction $\sfrac{1}{3}$, it should be done in terms of this point.
\end{quote}
The use of the word representative in this statement is remarkable, because representation usually works the other way around: $\sfrac{1}{3}$ represents a point on the number line, with $\sfrac{2}{6}$ 
representing the same point.

\subsubsection{Fractions as values and the BOCARDO syllogism}
Suppose one assumes that ``all fractions have a numerator and a denominator'', and that ``all fractions are
rational numbers'', then the BOCARDO syllogism yields that ``some rational numbers have a numerator and a denominator''. But rational numbers don't come with a numerator and a denominator, and therefore at least one of the assumptions must be invalid.\footnote{%
This point is made quite explicitly (in Dutch) in Goddijn~\cite{Goddijn2005} (in Dutch):
\begin{quote}
De rationale getallen zijn een geconstrueerd systeem van `invulantwoorden' op eerder onmogelijke 
vermenigvuldigopgaven. Zo'n antwoord `heeft', geen teller of noemer.
\end{quote} }

This argument is at first sight conclusive although a way out is to make the  assumption that a
after all a  (positive) rational $q$ can be issued a numerator $n$ and a denominator $d$, 
both positive naturals, with greatest common divisor equal to $1$ so that that $q = n/d$. I will refer to this choice of 
a numerator and a denominator of a rational number as the standard decomposition of a rational number.

A more principled approach is to abandon Aristotelian logic, or at least to develop an explanation of fractions in
which some applications of  the BOCARDO syllogism are not implied. Reiser~\cite{Reiser1936} provides clear indications that such alternatives may be needed. Details of non-Aristotelian logic, including the various syllogisms can be found in Vasil'ev~\cite{Vasilev2003}. Within a non-Aristotelian paradigm Geeraerts~\cite{Geeraerts1983} provides a case study of prototype theory and diachronic semantics. This work might wel be be applicable to the case of fractions: fraction as a prototype with one or more clusters of temporally interrelated meanings. An alternative step towards dismantling undesired consequences of syllogistic logic 
consists of contemplating a more sophisticated model of dialogue and interaction where it is not taken for 
granted that both agents in a 
discourse assign the same meaning to the same assertions, see for instance Kent~\cite{Kent1983}.

\subsubsection{The Simplifiable Fraction Anomaly}
\label{SimplFA}
Assuming that fractions are rational numbers creates via the BOCARDO syllogism the unexpected implication that some rationals are decomposable. Though perhaps unnatural this implication is by no means inconsistent. How to decompose rational numbers into relatively prime components has been indicated above.
 Now consider the class of simplifiable fractions (e.g. $2/4$ or $3/6$ etc.). I will assume that all simplifiable fractions are fractions. 
 Upon assuming that all fractions are rational numbers, 
one finds that all simplifiable fractions are rational numbers. BOCARDO now yields that 
some rational numbers are simplifiable. The latter assertion is wrong when rationals are decomposed in the standard manner because the standard decomposition yields a simplified fraction. As there is no plausible 
other way of assigning a numerator and a denominator to a rational number an anomaly appears from the mere assumption that fractions are rational  numbers in combination with the assumption that fractions are decomposable and the assumption that some fractions are simplifiable. I will refer to this anomaly as the Simplifiable Fraction Anomaly.

These three assumptions are met by the exposition on fractions of Van de Craats~\cite{Craats2010} who
writes that a fraction is the result (``uitkomst'') of a division of two integers. 
This wording amounts to the assertion that a fraction is a rational number. More specifically he indicates that $4:7 = \frac{4}{7}$ which 
suggests that $4:7$ is an expression which is evaluated (rewritten) into a value (result) $\frac{4}{7}$. 
As a consequence, the exposition in~\cite{Craats2010} features the Simplifiable Fraction Anomaly.

\subsubsection{Fraction naming}
When adopting the fractions as values paradigm names of 
numbers, and names of fractions, are likely to play a prominent role in fact to play the role played of what is termed a fraction in other approaches.  For instance in the fractions as values paradigm the name of a fraction, rather than the fraction itself, is supposed to be equipped with a numerator and a denominator.  

Rollnik argues at length that many approaches to fractions found in the literature and in existing teaching materials  lead to mistakes, imprecision, and even contradictions, and he argues that fractions as a values is the better choice. Similar but less comprehensive criticism is formulated by Opmeer in~\cite{Opmeer2005}. The proposal of Rollnik depends, however, on the ability to provide a useful account of names of rational numbers. 

Now unfortunately the very notion of a name is not so simple and its analysis has lead to intricate philosophical ramifications. 
For instance: is ``$\frac{1}{2}$ is smaller than $1$'' a so-called bare use of the 
name $\frac{1}{2}$ or is it a predicative use?
In~\cite{Englebretsen1981} one finds the position that the distinction between bare use and predicative use 
is a context-sensitive one, thus 
opposing a famous view by Burge in~\cite{Burge1971}, who made the suggestion that names are predicates.
It is possible to view $\frac{1}{2}$ as a predicate on descriptions of numbers?

  For recent work on names see Gray~\cite{Gray2012}. Gray explains that a philosophical theory of names seeks to explain what a formalised logical theory of names in the tradition of analytical philosophy intends to avoid: 
  the complications arising with the use of names in natural language. Following Gray's view in that matter the fractions as terms paradigm which I prefer, 
  intentionally, or at least consciously, avoids giving an account of naming rational numbers.

\subsection{Fractions as fracpairs in the literature}
The term fracpair is not used in the papers mentioned in this Paragraph, what matters here is that a fraction is identified with a pair of integers.

This definition of fractions as pairs with addition rule MFAR, though in the presence of zero and of negative numbers, and with a slightly different notation, is given by Gregorczyk in~\cite{Gregorczyk1974}. Kleene~\cite{Kleene2002} 
also views fractions as pairs of integers with the same rewrite rules for operations. 
Kleene explicitly states that fraction equality requires both components of the pair to be the same.
 In Ni~\cite{Ni2001} it is stated that
\begin{quote}
Every fraction belongs to a single equivalence class generated by a single multiplicative equation, and each equivalence class defines a distinct rational number.
\end{quote}

The course material written by  Notten \& Verheij~\cite{NottenV2013} is compatible with a perspective 
where the fraction base type is fracpair, and the fracsign target type is fracpair as well.

In Prediger~\cite{Prediger2008} in Table 1 it is mentioned that among the 
necessary change of conceptions on finds that
\begin{quote}
Existence of many fractions representing the same fractional number
\end{quote} whereas for natural numbers: she speaks of:
\begin{quote}
Unique relation between number and symbolic representation
\end{quote}
Clearly $1/2$ and $2/4$ are considered different fraction, while the 
paper does not indicate whether or not $(1+2)/3$ and $3/3$ are considered different fractions. 

In Behr \& Post~\cite{BehrP1992} fraction equivalence is a central notion. 
In~\cite{BehrP1992} the phrase general fraction is used
in connection with an arbitrary pair of integers. From these facts I conclude that Behr views fractions as pairs.
This is not very obvious from the text, however, as the basic fraction concept, which also plays a central role in the paper, seems itself to be a spectrum of more precise subconcepts.

Adjiage \& Pluvinage~\cite{AdjiageP2007} use the phrase ``fractional writing'' for what has been termed a fracpair above (while taking zero and negative numbers into account). Following Duval (see e.g.~\cite{RDuval2000}) they 
refer to the variety of semiotic registers for fractions with fractional writing constituting one semiotic register and 
decimal writing and linear scale constituting different semiotic registers made us of in the same paper.

%

%
%
%
 A remarkably detailed 
description of fractions as fracpairs occurs in~Van Engen~\cite{Engen1960}. 
Van Engen views fractions as pairs,
to be distinguished from so-called rate pairs, and the distinction depends on the operators which can be applied. 
Rate pairs are not equipped with arithmetical operations, 
while fractions are. 
Van Engen's paper~\cite{Engen1960} is the most closely related work to the current paper I have found in the literature on fractions.
Glennon \& Callahan comment on  Van Engen's proposal in~\cite{GlennonC1986} (p93) 
stating that it considers fractions to be numbers, a comment which I do not agree with.

\section{More notational issues}
Fractions as a subject is quite extended, and leads to many notational issues, some of which I will discuss in this Section.
\subsection{How to write the division operator}
\label{DivOps}
In~\cite{Craats2010} Van de Craats states that with a fraction, say $\displaystyle \frac{3}{7}$, three calculations (``sommen'' in Dutch) are 
to be associated: $\displaystyle 4:7 = \frac{4}{7}, 7 \times   \frac{4}{7} = 4$ and $\displaystyle 4 \times \frac{1}{7}=   \frac{4}{7}$. It appears that the horizontal bar does represent a notational connective rather than a function symbol and that in~\cite{Craats2010} only division written with infix symbol ``$:$'' and multiplication, written ``$\times$'' occur as function names. However, later
in the same text one finds $\displaystyle \frac{5}{7} \times \frac{3}{4}=  \frac{5 \times 3}{7 \times 4}$ an identity which
seems to indicate that the horizontal bar denotes a function symbol after all. If that were the case, however, 
the distinction between $4:7$ and $\displaystyle  \frac{4}{7}$ becomes even more puzzling. Are $4:7$ and 
$\displaystyle  \frac{4}{7}$ synomyms in the same way as $\displaystyle  \frac{4}{7}$ and 
$4/7$ are merely different ways of writing the same (according to~\cite{Craats2010})?
In any case, as far as I can see  the Simplifiable Fraction Anomaly  
can be reconstructed in the setting of~\cite{Craats2010}.

\subsection{Integer division on positive natural numbers}
\label{DoPNN}
When working with whole numbers students are supposed to solve exercises like ``$ 10 : 2 = \ldots$''. 
The right answer is supposedly to write $5$ on the dotted part of the text fragment (``incomplete-equation-sign''). 
A further question may read ``$ 11 : 2 = \ldots \mathrm{rem} \ldots$'', and the correct answer is supposedly to write
5 and 1 on both dotted pars of the sign (in that order, when reading the sign from left to right).

In the second case the equality sign seems not to express the equality of value between two expressions, 
one on its left-hand side and one on its right-hand side. The view that $ 11 : 2 = \ldots \mathrm{rem} \ldots$ 
constitutes an 
equality may be maintained if the outcome of evaluating $ 11 : 2 $ is taken to be the pair $(5,1)$. 
Doing so, however, 
has the somewhat 
unfortunate consequence that for instance $2:3 = 3:4$. 
Because the latter identity will be denied once ratios are discussed another reading of 
$ 11 : 2 = 5 ~\mathrm{rem} ~1$ may be both preferable and in fact also more plausible: 
a four place relation on (positive) natural numbers written by means of a mixfix notation.  
In fact this is the only interpretation of  $ 11 : 2 = 5 ~\mathrm{rem}~ 1$
 which makes sense to me, as writing say $ (11 : 2 = 5 ~\mathrm{rem}) ~1$, or  writing 
  $ 11 : 2 =_{rem} 5 ~\mathrm{rem} ~1$, or writing $ 11 : 2 =_{rem~1} 5~\mathrm{rem} ~1$, or writing 
  $ 11 : 2 =_{rem~1} 5$ is not of much help, given the fact that in each case a non-transitive form of equality results.
  
 It might be preferable to write $10:2 => 5$ in advance of writing $10:2 = 5$ and to write 
 $ 11 : 2 => 5 ~\mathrm{rem}~ 1$ instead of $ 11 : 2 = 5 ~\mathrm{rem}~ 1$ thus
 avoiding the use of the equality sign in a way which on the long run may cause confusion. 
 
$10:2 =>5$ is read as ``the processing of $5:2$ is successful and produces $2$'', and 
$ 11 : 2 => 5 ~\mathrm{rem}~ 1$ is read as ``the processing of $5:2$ is successful and produces $2$ with a remainder equal to $1$''.
 
Now it may be stated that  if for some $n,m$ and $k$ it is known that $n:m =>k$, 
the equation $n:m = k$ holds as well, and moreover one may say the there is no (positive natural) number $k$ 
such that $11:2 =>k$.
 
Using $\backslash$ as an infix operator symbol for integer division and using $\%$ as an infix operator symbol 
for remainder one finds:
$n:m = > n \backslash m~\mathrm{rem}~n\,\%\,m$. Moreover $11 \backslash 2 = 5$ and $11\, \%\, 2 = 1$. 
Further 
$10\backslash 2$, $10\,\%\, 2$, $11 \backslash 2 $, and $11\, \%\, 2 $  may be considered expressions because it is 
entirely clear how these can be evaluated, this in contrast with $10:2$ which allows successful 
evaluation to a single value  ``by accident'' only, unlike and $11:2$
which defeats evaluation to any value in a known domain of objects (unless of course pairs count as values).

%

\subsection{Decimal notation and percentages}
Decimal fractions may be understood as interpretations of decfracsigns, that is signs in the range of syntax 
known for decimal fractions. The collection decfracsigns is disjoint from fracsigns as each fracsign is supposed to have a division operator as its leading operator symbol. Casting a decfracsign to a fracvalue is trivial. Casting it to  
a fracpair (written as a fracpairsign, that is a fracsign for a simple fraction) or a fractriple (written as a mixfracsign for a mixed fraction) is harder. I will discuss some examples only. 
Reading the decfracsign $0.8$ as a fracpair yields $\displaystyle \frac{8}{10}$
rather than $\displaystyle \frac{4}{5}$. Reading $0.85$ as a fracpair yields $\displaystyle \frac{85}{100}$
rather than $\displaystyle \frac{17}{20}$, and reading $2.4$ as a fractriple yields 
$\displaystyle 2\frac{4}{10}$ rather than $\displaystyle 2\frac{2}{5}$.

Each fraction content policies may be adapted to a systematic inclusion of decimal notation. 
There will be several options for doing so and just as with the case of mixed fractions I hold that this topic is probably 
best viewed in terms of fraction as constituting a dynamic notion.

Working with percentages may be supported by casting say $50 \%$ to 
$\displaystyle \frac{50}{100}$ rather than say $\displaystyle \frac{1}{2}$ and casting 
$\displaystyle 7\frac{1}{2} \%$ as well as $\displaystyle 7.5 \%$ to a fracpair as $\displaystyle \frac{75}{1000}$.
Finally 
$\displaystyle 125 \%$ may be cast to the fracpair $\displaystyle \frac{125}{100}$ and to the fractriple 
$\displaystyle 1\frac{25}{100}$. Casting a percentage to a fracterm is always trivial, however,
for instance $\displaystyle 125\frac{2}{4} \%$ may be first cast to the fracterm
$\displaystyle \frac{125\frac{2}{4}}{100}$ and subsequent casting may be performed from there.

Yet another view is first to cast a percentage to a decimal notation and only thereafter to cast it to 
a fracpair, a fractriple, or a fracterm.

Most sources explaining how to convert percentages to fractions will transform to a simplified 
simple fraction or to a corresponding mixed fraction, rather than to return a fraction which, as a fraction, is most closely resembling the given percentage.

\subsection{Multiset ratio model with characters}
I will first consider multisets of alphabetical characters and character strings, for instance $V=\{a,a,a,a,A,A,A,A,A,A\}$. 
$V$ might be depicted as a collection of 4 green circles with 6 yellow circles, or in any other suitable manner. Now 
the following assertions about $V$ and its elements may be used to link $V$ with the terminology of ratios. 
I will assume that ratios are understood as fracpairs. 
\begin{enumerate}
\item The ratio of lower case characters to upper case characters is $4/6$. 
\item Alternative formulations of the same state of affairs: 
the ratio of lower case characters to upper case characters: 
(i) is given by $4/6$, 
(ii) corresponds to $4/6$, 
(iii) matches with $4/6$, 
(iv)  is equivalent to $2/3$, or 
(v) is equivalent to $6/12$.

\item The ratio of small characters to capital characters \emph{matches best} with  $4/6$, the match is better than with
the match with $2/3$ and the match with $8/12$.
\item The ratio of lower case characters to all of $V$ is $4/10$.
\item The ratio of lower case characters to all of $V$ is best given by $4/10$.
\end{enumerate}

These assertions provide a way of combining terminology for the construct fractions of numbers with terminology for 
ratios in connection with multisets. The multiset approach belongs to the part-whole fraction construct.
 Novillis~\cite{Novillis1976} provides a survey of fraction perspectives linked with the part-whole fraction construct. 
On p143 Novillis writes that:
\begin{quote}
Another inference that can be made is that students do not come into contact with an adequate number of negative instances of the fraction concept.
\end{quote} which indicates the supposed existence of a boundary between fraction and non-fraction which might 
open up a new dimension of exploration altogether. What are meaningful non-fractions?

\subsubsection{Multiset union and ratio union}
\label{MUandRU}
Ratios allow a union operator which is not compatible with equivalence, that is the union of pairwise 
equivalent ratios need not produce a pair of equivalent ratios.
\begin{definition}Ratio union operator: $\displaystyle \frac{a}{b} \cup\frac{c}{d} = \frac{a+c}{b+d}$.
\end{definition}
Ratio union describes the best ratio of the union of two multisets each consisting of two kinds of elements as a function of the best ratio for its components.
Technically ratio union is just vector addition, taking to the extreme the point of view that a fracpair is just  a vector of dimension 2. This operation has been investigated in detail in Mochon~\cite{Mochon1993}. 
Mochon 
analyses cases where it makes sense to combine fractions by way of vector addition. He writes 
$\frac{a}{b}$``+'' $\frac{c}{d}$ instead of $\frac{a}{b} \cup\frac{c}{d}$. For instance let 
$W = \{a,d,d,B,B,C,D,E\}$ then the best ratio of lower case to upper case elements for $W$  is $3/5$ and 
so the best ratio for $V \cup W$ is $7/11$.

\subsubsection{Multisets of coloured tokens}
Let $C=\{g,b,r,y,p,w\} $ serve as a collection of constants for tokens with various colours. 
I now consider multisets over $C$ on which the following operations are defined: $\#(V)$ is the cardinality of $V$, 
For as set $U \subseteq V$, $\#_U(V)$ counts the number of elements in $V$ which are also in $U$. Thus
$\#_C(V) =\#(V)$. $\mathit{ratio}_U(V) = \#_U(V):\#(V)$.
Further like ratios can be added by means of a version of the quasicardinality rule. I will use the operator symbol
most popular for ratios in this case:
\begin{definition} Addition of like ratios: $(a: b) + (c:b) = (a+c):b$.
\end{definition}
\noindent Now for different colours, say $g$ and $b$, the following equation expresses a count for two different colours in terms of the corresponding counts for a single colour: $\#_{\{g,b\}}(V) = \#_{\{g\}}(V) +\#_{\{b\}}(V) $. From this we find:
$\mathit{ratio}_{\{g,b\}}(V) = \mathit{ratio}_{\{g\}}(V) + \mathit{ratio}_{\{b\}}(V)$. This identity, when explained in words, provides a convincing application of the quasicardinality rule for ratios.

\subsubsection{Multisets in the literature on fractions.} Thompson \& Lambdin~\cite{ThompsonL1994} consider graphical displays consisting of sets of circles with different colours. 
It is indicated that a diversity of interpretations of such displays in the content of fractions and part-whole relations exist 
and~\cite{ThompsonL1994} states that rather than requiring a student to choose a unique ``correct'' way of reading such pictures 
students must be made aware of the full diversity of ways of reading such pictures. Said diversity  also includes a
reading which corresponds to $\sfrac{2}{5}+\sfrac{3}{4} = \sfrac{6}{9}$ where $+$ stands for the combination 
of collections of coloured items: 2 green items out of 5 items combined with 3 green items out of 4 items is 
6 green items out of 9.
Gould~\cite{Gould2013} also surveys different interpretations of diagrams (coloured pizza parts), 
but concludes that many of these are misleading, or represent student misunderstanding.

\subsection{On different decimal fractions}
\label{DecNot}
The difficulty concerning the concept of fractions comes about 
in the hesitation felt by many authors to portray say $2/3$ and $4/6$ as different fractions and to develop a fraction talk which takes this difference into account in a systematic and explicit manner.
The issue is not specific for ``fractional notation'', however. In decimal notation one may wonder in what sense $7.5$ and $7.50$ are different. Clearly both notations for a rational feature different numbers of digits. This comes with the connotation that, if say  rational number $r$ is specified by the equation $r = 7.50$ while $s$ is specified by $s=7.5$, the specification of $r$ may be understood  
to provide higher absolute precision than the specification of $s$ is providing (viz. $0.005$ rather than $0.05$). 

Translating conventions from decimal notation to fractional notation  one might assume that $2/3$ features an implicit
absolute precision of $1/6$ while $4/6$ features an implicit absolute precision of $1/12$ but such conventions seem to be absent in the literature on fractions. Taking such (potential) conventions unconventionally seriously the decimal to fraction conversion of $7.5$ yields $75/1000$ while the conversion of $7.50$ would yield $750/10000$. To the best of my knowledge it is convention to simplify fractions after decimal to fraction conversion with the effect that both conversions yield $3/4$. It appears that ``having a different number of decimals'' separates decimal notations in 
a manner for which no counterpart is known or used for ordinary fractional notation.

\section{Concluding remarks on defining fractions}
 The very idea of a principled syntax dismissive attitude is foreign to  informatics. 
Nevertheless a 
syntax dismissive (and naming dismissive) approach seems indispensable for the introduction of 
the most basic ideas in mathematics. Once these basic ideas and structures, including booleans, positive naturals, 
naturals, and integers, have been made available, possibly in a partially syntax permissive paradigm, 
further development may increasingly take place in a syntax permissive setting. Fractions and rational numbers 
offer a ``battleground'' where naming permissiveness and syntax permissiveness may enter traditionally
mathematical territory. 

A contrast between informatics and mathematics is present in diverging attitudes towards expressions. Whereas in mathematics an expression comes with a meaning and the determination of its value is considered a critical task, in informatics 
the determination of its meaning, in advance of the 
computation of a particular value, is considered a meaningful challenge.
For instance $1\div 0$ qualifies as an expression for a computer scientist, the meaning of which may yet be determined. 
For a mathematician the very fact that the meaning of $1\div 0$ is in doubt disqualifies it as an expression. 
Viewing $1/2$ as an expression suggests that it can be evaluated and that it has a value.
But as it is natural to perceive $1/2$ as a value directly, is classification as an expression loses credibility.
With $2 + 3$ no such problem arises assuming that one takes $5$ as a value for it.

\subsection{Defining fractions compared with other defininitions}
 Suppose one is interested in writing a survey of human transportation. 
 A text is written about cars, bicycles, trains, airplanes, ships, and busses. Now the question is posed to introduce \emph{means of human transportation} (MoHT), because this notion if felt to constitute a useful generalisation of specialised forms of transportation. 
 However, when introducing  MoHT in a rigorous manner new questions arise. For instance: is an elevator is an MoHT, is an escalator an MoHT,  is a parachute an MoTH, and is a horse an MoHT? 
  
One may write convincingly about human transportation and yet be unable or unwilling to commit to any rigorous definition of a means of human transportation.  In a comparable manner it seems to be an illusion to infer the existence of a rigorously definable concept of fraction from either the importance of fractions or from the large volume of existing writings about fractions.

\subsection{Fractions defeating definition in advance of use?}
Answering the question ``what is a word'' by making use of written natural language is difficult. 
In a proposed definition words will be used to define words,  and a circularity is bound to arise. 
The notion of a word is so close to natural language itself that expecting to
read about words in advance of learning to use words is implausible if not unreasonable. 
One can learn about bats as a concept biology in advance of searching the sky for such animals, 
but one cannot learn about words as linguistic concept in advance of dealing with words.

I hold that ``fraction'' shares with ``word'' the feature that defining it (i.e. providing an answer to the
question ``what is a fraction'') before using it is unfeasible: 
teaching fractions cannot plausibly start with providing a valid definition of a fraction in simpler terms. 
Therefore, only after having been equipped with an initial grasp of fractions and numbers, 
which is likely to feature gaps as well as inconsistencies,  one may proceed to a second stage of understanding
in which adequate definitions of number and fraction are leading.
%
%

%
%
 
\subsection{Acknowledgements}
The idea to formulate the simplifiable fraction anomaly and to investigate its role in an analysis of definitions of fractions came about as a consequence of an exchange of emails with Joost Hulshof (VU Univ. Amsterdam) 
about the merits of the definitions in teaching material on fractions in Van de Craats~\cite{Craats2010}. 
This exchange had been triggered by remarks 
made in various emails sent to both of us by Thomas Cool 
(independent researcher, The Netherlands). 
Unsurprisingly Joost and I came nowhere near an 
agreement on the issues at hand. 
Ben Wilbrink (independent researcher on educational psychology, the Netherlands) brought  Wu~\cite{Wu2014} to my attention.   
Jeroen van Wier 
(MSc Logic student at UvA) has suggested several textual improvements. Input from Anngha Nugraha
(visiting the UvA) led to the remarks on  decimal notation.
In different stages of the work  I had useful discussions and email exchanges concerning fractions and rational numbers with the following persons: Inge Bethke (UvA), Andrea Haker (HvA), Dimitri Hendriks (VUA), Kees Middelburg (uvA), Alban Ponse (UvA), Stefan Rollnik,  
Anja Sicking, John Tucker (Univ. of Wales Swansea), and Albert Visser (UU). 
These persons bear no responsibility for what I wrote.

\end{document}